\documentclass[11pt,reqno]{amsart}
\usepackage{graphicx}
\usepackage{verbatim}
\usepackage{textcomp}
\usepackage{amssymb}
\usepackage{cite}
\usepackage{amsmath}
\usepackage{latexsym}
\usepackage{amscd}
\usepackage{amsthm}
\usepackage{mathrsfs}
\usepackage{bm}
\usepackage{url}
\usepackage{hyperref}
\usepackage{bookmark}
\allowdisplaybreaks[3]
\newtheorem{thm}{Theorem}[section]
\newtheorem{corr}[thm]{Corollary}
\newtheorem{lem}[thm]{Lemma}

\theoremstyle{definition}

\theoremstyle{remark}
\newtheorem{rem}{Remark}[section]
\numberwithin{equation}{section}
\begin{document}
\title[Rigidity of Einstein metrics on complete manifolds]
{Rigidity of Einstein metrics as critical points of some quadratic curvature functionals on complete manifolds}

\author{Guangyue Huang}
\author{Yu Chen}
\author{Xingxiao Li}

\address{Department of Mathematics, Henan Normal
University, Xinxiang 453007, People's Republic of China}
\email{hgy@henannu.edu.cn (G. Huang)}
\email{yuchern@126.com (Y. Chen)}
\email{xxl@henannu.edu.cn (X. Li)}

\thanks{Research supported by NSFC (Nos. 11401179, 11371018, 11671121).}

\begin{abstract}
In this paper, we consider some rigidity results for the Einstein metrics as the critical points
of some known quadratic curvature functionals on complete manifolds, characterized by some point-wise inequalities. Moreover, we also provide rigidity results by the integral inequalities involving the Weyl curvature, the trace-less Ricci curvature and the Sobolev constant, accordingly.
\end{abstract}



\maketitle

\section{Introduction}

In this paper, we always assume that $M^n$ is a complete manifold of dimension $n\geq3$ and $g$ is a Riemannian metric on $M^n$ with the Riemannian curvature tensor $R_{ijkl}$, the Ricci tensor $R_{ij}$ and the scalar curvature $R$.
It is well-known that Einstein metrics are critical points for the Einstein-Hilbert functional
\begin{equation*}
\mathcal {H}=\int_{M}R
\end{equation*}
on the space of unit volume metrics $\mathscr{M}_1(M^n)$. In \cite{Catino2015},
Catino considered the following family of quadratic curvature functionals
\begin{equation}\label{Int-1}
\mathcal{F}_t\,=\,\int_M |R_{ij}|^{2}+t\int_M
R^{2},\quad t\in \mathbb{R}
\end{equation}
which are also defined on $\mathscr{M}_1(M^n)$, and proved some related rigidity results. Furthermore,
it has been observed in \cite{Bess2008} that every
Einstein metric is a critical point of $\mathcal{F}_t$ for all $t\in \mathbb{R}$ (see \eqref{2-lemmas-16} in Section 2). But the converse of this conclusion is not true in general.

Therefore, it is natural to study canonical metrics which arise as solutions of Euler-Lagrange equations for more general curvature functionals. For instance, Anderson \cite{Anderson2001} proved that
every complete three-dimensional critical metric for the Ricci functional $\mathcal{F}_0$ with non-negative
scalar curvature is flat. In \cite{Catino-PorcAMS2016}, Catino gave a characterization of complete
critical metrics for the functional
\begin{equation*}
\mathcal {H}=\int_{M}R^2
\end{equation*}
with non-negative scalar curvature in every dimension. For the $\sigma_2$-curvature functional on $M^3$, Catino \cite{Catino-Pacific2016} proved that flat metrics are the only complete metrics with non-negative scalar curvature.

By imposing some inequalities conditions involving the Weyl tensor $W$ and the traceless Ricci tensor $\mathring{\rm Ric}$, with components $W_{ijkl}$ and $\mathring{R}_{ij}$, respectively, we are able to prove a few rigidity results for the Einstein metrics considered as the critical points of the functional $\mathcal{F}_t$ given by \eqref{Int-1}. 

Firstly, we define two constants $D_n$ and $E_n$ by
\begin{align}D_n=&\begin{cases}\frac{4}{n^2},\ \ \ \ \ \quad \quad {\rm if}\ n=3,4,5,6;\\
\Big(\frac{4}{n}-1\Big)^2,\ \  {\rm if}\ n\geq7;
\end{cases}\label{dn}\\
E_n=&\begin{cases}\sqrt{6},&\text{if\ }n=4;\\
\frac{4(n-1)}{n(n-2)}\Big(\frac{n-2}{\sqrt{n(n-1)}} +\frac{n^2-n-4}{2\sqrt{(n-2)(n-1)n(n+1)}}\Big)^{-1},&\text{if\ }n\geq 6.
\end{cases}\label{En}
\end{align}

Then our main results in this paper can be stated as follows:

\begin{thm}\label{thm1-1}
Let $(M^n,g)$ be a complete Riemannian manifold of dimension $n\geq 3$ with
positive scalar curvature and
\begin{align}\label{111thm-Int-1}
\int_M|\mathring{\rm Ric}|^2<\infty.
\end{align}
Suppose that $g$ is a critical metric for
the functional $\mathcal{F}_t$ over $\mathscr{M}_1(M^n)$ with
\begin{equation}\label{1111thm-Int-1}
\begin{cases}t\leq-\frac{5}{12},\ \ \ \ \quad \quad \text{when}\ n=3;\\
t\leq-\frac{13n-1}{50(n-1)},\  \quad \text{when}\ n\geq4. \\
\end{cases}
\end{equation}
If
\begin{align}\label{1thm-Int-4}
\Big|W&+\frac{\sqrt{n D_n}}{\sqrt{2}(n-2)}\mathring{{\rm Ric}} \mathbin{\bigcirc\mkern-15mu\wedge} g\Big|\leq-\frac{2\sqrt{2}[(n-2)+n(n-1)t]}{n\sqrt{(n-1)(n-2)}}R,
\end{align}
then $M^n$ is Einstein. In particular, when $n=3$, $M^3$ must be of constant positive sectional curvature; when $n\geq4$,
if the parameter $t$ also satisfies that
\begin{equation}\label{1111thm-Int-6}
\begin{cases}-\frac{\sqrt{2}+1}{6\sqrt{2}}<t\leq-\frac{51}{150},\ \ \ \, \quad \quad\quad \ \quad \quad \ \quad \ \ \quad \quad \text{when}\ n=4;\\
-\frac{3}{20}\Big(1+\frac{\sqrt{15}}{8}\Big)<t\leq-\frac{8}{25},\ \ \ \ \quad \quad\quad \quad \ \quad\quad \text{when}\ n=5; \\
-\Big(\frac{n-2}{n(n-1)}+\frac{1}{2\sqrt{2}nC_n}\sqrt{\frac{n-2}{n-1}}\Big)<t\leq-\frac{13n-1}{50(n-1)},\   \text{when}\ n\geq6,
\end{cases}
\end{equation}
then $M^n$ must be of constant positive sectional curvature.

\end{thm}

We recall that the Sobolev constant $Q_g(M)$ is defined by
\begin{equation}\label{1thm-Int-6}
\aligned
Q_g(M)=\inf\limits_{0\neq u\in C_0^{\infty}(M)}\frac{\int_M\Big(|\nabla u|^2+\frac{n-2}{4(n-1)}R_gu^2\Big)}{(\int_M|u|^{\frac{2n}{n-2}})^{\frac{n-2}{n}}}.
\endaligned\end{equation}
For example, any simply connected complete manifold with $W\equiv 0$ has positive Sobolev constant (see \cite{Schoen-yau1988}).
Moreover, it is easy to see from \eqref{1thm-Int-6} that
\begin{equation}\label{1thm-Int-7}
\aligned
Q_g(M)\Big(\int_M|u|^{\frac{2n}{n-2}}\Big)^{\frac{n-2}{n}}\leq\int_M\Big(|\nabla u|^2+\frac{n-2}{4(n-1)}R_gu^2\Big)
\endaligned\end{equation}
for any smooth $L^{\frac{2n}{n-2}}$-function $u$. With the help of \eqref{1thm-Int-7},
we can prove the following rigidity result:

\begin{thm}\label{thm1-2}
Let $(M^n,g)$ be a complete Riemannian manifold of dimension $n\geq 3$ with
positive scalar curvature and $g$ be a critical metric for
the functional $\mathcal{F}_t$, $t<-\frac{1}{2}$, over $\mathscr{M}_1(M^n)$ satisfying
\begin{align}\label{222thm-Int-1}
\int_M|\mathring{\rm Ric}|^2<\infty.
\end{align}
If $Q_g(M)>0$ and
\begin{equation}\label{222thm-Int-3}
\Big(\int_M\Big|W+\frac{\sqrt{2}}{\sqrt{n}(n-2)}\mathring{{\rm Ric}} \mathbin{\bigcirc\mkern-15mu\wedge} g\Big|^{\frac{n}{2}}\Big)^{\frac{2}{n}}<\sqrt{\frac{n-1}{2(n-2)}}Q_g(M),
\end{equation}
then $M^n$ is Einstein. In particular,

(1) when $n=3,4,5$, $M^n$ must be of constant positive sectional curvature;

(2) when $n\geq6$, if we replace \eqref{222thm-Int-3} with
\begin{equation}\label{222thm-Int-4}
\Big(\int_M\Big|W+\frac{\sqrt{2}}{\sqrt{n}(n-2)}\mathring{{\rm Ric}} \mathbin{\bigcirc\mkern-15mu\wedge} g\Big|^{\frac{n}{2}}\Big)^{\frac{2}{n}}<E_nQ_g(M),
\end{equation}
then $M^n$ must be of constant positive sectional curvature.
\end{thm}

Since there is no complete noncompact Einstein manifold with positive constant scalar curvature, the following results follows easily from the above two theorems:

\begin{corr}\label{4corr1-2}
Suppose that $(M^n,g)$ is a complete noncompact Riemannian manifold of dimension $n\geq 3$ with
positive scalar curvature and that $g$ is a critical metric for
the functional $\mathcal{F}_t$ over $\mathscr{M}_1(M^n)$ with $t$ satisfying \eqref{1111thm-Int-1}.
If \eqref{1thm-Int-4} holds, then it must hold that
\begin{align}\label{corr-Int-1}
\int_M|\mathring{\rm Ric}|^2=\infty.
\end{align}

\end{corr}

\begin{corr}\label{4corr1-1}
Suppose that $(M^n,g)$ is a complete noncompact Riemannian manifold of dimension $n\geq 3$ with
positive scalar curvature and that $g$ is a critical metric for
the functional $\mathcal{F}_t$ over $\mathscr{M}_1(M^n)$ with $t<-\frac{1}{2}$.
If $Q_g(M)>0$ and \eqref{222thm-Int-3} holds, then it must hold that
\begin{align}
\int_M|\mathring{\rm Ric}|^2=\infty.
\end{align}

\end{corr}

\begin{rem}\label{1Rem-2}
Noticing that if $M^n$ has constant sectional curvature, then both the Weyl tensor $W$ and the traceless Ricci tensor $\mathring{\rm Ric}$ are identically equal to zero, implying that the function
\begin{align}\label{Rem-1}
\Big|W+\frac{\sqrt{n D_n}}{\sqrt{2}(n-2)}\mathring{{\rm Ric}} \mathbin{\bigcirc\mkern-15mu\wedge} g\Big|
\end{align}
is identically zero. But in general the right hand side of \eqref{1thm-Int-4} is strictly positive by \eqref{1111thm-Int-1} and the assumption that $R>0$. Therefore, Theorem \ref{thm1-1} actually reveals a pinching phenomenon for the function \eqref{Rem-1}. Similarly, when $n=3,4,5$, Theorem \ref{thm1-2} also discloses a different pinching phenomenon for the quantity
\begin{align}\label{Rem-2}
\int_M\Big|W+\frac{\sqrt{2}}{\sqrt{n}(n-2)}\mathring{{\rm Ric}} \mathbin{\bigcirc\mkern-15mu\wedge} g\Big|^{\frac{n}{2}}.
\end{align}

\end{rem}

\begin{rem}\label{1Rem-1}
The first two authors of the present paper have previously proved in \cite{MHLC2018} some rigidity results on Einstein metrics as the critical points
of the same quadratic curvature functionals for the case of compact manifolds. Our theorems are then the generalizations of those conclusions to the case of complete manifolds.
\end{rem}

\section{Some necessary lemmas}

Recall that the Weyl curvature $W_{ijkl}$ of a Riemannian manifold $(M^n,g)$ with $n\geq3$ is related to the Riemannian curvature $R_{ijkl}$ by
\begin{align}\label{2-lemmas-1}
W_{ijkl}=&R_{ijkl}-\frac{1}{n-2}(R_{ik}g_{jl}-R_{il}g_{jk}
+R_{jl}g_{ik}-R_{jk}g_{il})\notag\\
&+\frac{R}{(n-1)(n-2)}(g_{ik}g_{jl}-g_{il}g_{jk}).
\end{align}
Since the traceless Ricci curvature $\mathring{R}_{ij}=R_{ij}-\frac{R}{n}g_{ij}$,
\eqref{2-lemmas-1} can be written as
\begin{align}\label{2-lemmas-2}
W_{ijkl}=&R_{ijkl}-\frac{1}{n-2}(\mathring{R}_{ik}g_{jl}-\mathring{R}_{il}g_{jk}
+\mathring{R}_{jl}g_{ik}-\mathring{R}_{jk}g_{il})\notag\\
&-\frac{R}{n(n-1)}(g_{ik}g_{jl}-g_{il}g_{jk}).
\end{align}
Furthermore, the Cotton tensor $C$ is defined to be with the components
\begin{align}\label{2-lemmas-3}
C_{ijk}:=&R_{kj,i}-R_{ki,j}-\frac{1}{2(n-1)} (R_{,i}g_{jk}-R_{,j}g_{ik})\notag\\
=&\mathring{R}_{kj,i}-\mathring{R}_{ki,j} +\frac{n-2}{2n(n-1)}(R_{,i}g_{jk}-R_{,j}g_{ik}),
\end{align}
where the indices after a comma denote the covariant derivatives.

We define the $(a,b)$-Cotton tensor $C^{(a,b)}$ to be with the components as follows:
\begin{equation}\label{2-lemmas-4}\aligned
C_{ijk}^{(a,b)}=\mathring{R}_{kj,i}-a\mathring{R}_{ki,j}+\frac{n-2}{2n(n-1)}(R_{,i}g_{jk}-bR_{,j}g_{ik}),
\endaligned
\end{equation}
where $a,b$ are two real constants. Clearly, the $(1,1)$-Cotton tensor $C_{ijk}^{(1,1)}$ is exactly the Cotton tensor $C_{ijk}$ defined by \eqref{2-lemmas-3}. In order to state our results, we first introduce the following cut-off function.
We denote by $p\in M^n$ and $B_r$ a fixed point and the geodesic sphere of $M^n$ of radius $r$ centered at $p$, respectively.  Let $\phi_r$ be the nonnegative cut-off function defined on $M^n$ satisfying
\begin{equation}\label{2-lemmas-5}
\phi_r=
\begin{cases}1,\ \ \quad {\rm on}\ B_r\\
0,\ \ \quad {\rm on}\ M^n\backslash B_{r+1}
\end{cases}
\end{equation}
with $|\nabla \phi_r|\leq2$ on $B_{r+1}\backslash B_{r}$.

First, we give an integral estimate on complete manifolds as follows:

\begin{lem}\label{Lemma22}
Let $M^n$ be a complete manifold. Then for any two positive constants $\epsilon_1, \epsilon_2$, we have
\begin{align}\label{2-lemmas-6}
&\int_M|\nabla \mathring{\rm Ric}|^2\phi_{r}^2\notag\\
\geq&\frac{1}{a^2+1+\epsilon_1}\int_M\Bigg[2a\Big(
W_{ijkl}\mathring{R}_{jl}\mathring{R}_{ik}
-\frac{n}{n-2}\mathring{R}_{ij}\mathring{R}_{jk}
\mathring{R}_{ki}
-\frac{1}{n-1}R|\mathring{\rm Ric}|^2\Big)\notag\\
&+|C^{(a,b)}|^2+\Big(-\frac{(n-2)^2(b-1)^2}{4n(n-1)^2}+\frac{(n-2)^2a}{2n(n-1)}-\epsilon_2\Big)|\nabla R|^2\Bigg]\phi_{r}^2\notag\\
&-\Big(\frac{4}{\epsilon_1}+\frac{(n-2)^2}{n\epsilon_2}\Big)\frac{a^2}{a^2+1+\epsilon_1}\int_M|\mathring{\rm Ric}|^2|\nabla \phi_{r}|^2.
\end{align}

\end{lem}

\proof From \eqref{2-lemmas-4}, we have
\begin{align}\label{2-lemmas-7}
|C^{(a,b)}|^2=&\sum_{i,j,k}|\mathring{R}_{kj,i}-a\mathring{R}_{ki,j}|^2
+\frac{(n-2)^2}{4n^2(n-1)^2}\sum_{i,j,k}|R_{,i}g_{jk}-bR_{,j}g_{ik}|^2\notag\\
&+\frac{n-2}{n(n-1)}(\mathring{R}_{kj,i}-a\mathring{R}_{ki,j})(R_{,i}g_{jk}-bR_{,j}g_{ik})\notag\\
=&(1+a^2)\sum_{i,j,k}|\mathring{R}_{ij,k}|^2-2a\mathring{R}_{kj,i}\mathring{R}_{ki,j}
+\frac{(n-2)^2}{4n^2(n-1)^2}[n(1+b^2)\notag\\
&-2b]|\nabla R|^2-\frac{(n-2)^2}{2n^2(n-1)}(a+b)|\nabla R|^2\notag\\
=&(1+a^2)\sum_{i,j,k}|\mathring{R}_{ij,k}|^2-2a\mathring{R}_{kj,i}\mathring{R}_{ki,j}\notag\\
&+\frac{(n-2)^2}{4(n-1)^2}\Big(\frac{(b-1)^2}{n}-\frac{2(n-1)a}{n^2}\Big)|\nabla
R|^2,
\end{align}
where we used the second Bianchi identity $\mathring{R}_{ij,j}=\frac{n-2}{2n}R_{,i}$.
Multiplying both sides of \eqref{2-lemmas-7} by $\phi_r^2$ and integrating it yield
\begin{align}\label{2-lemmas-8}
\int_M|C^{(a,b)}|^2\phi_r^2=&(1+a^2)\int_M|\nabla \mathring{\rm Ric}|^2\phi_r^2-2a\int_M\mathring{R}_{kj,i}\mathring{R}_{ki,j}\phi_r^2\notag\\
&+\frac{(n-2)^2}{4(n-1)^2}\Big(\frac{(b-1)^2}{n}-\frac{2(n-1)a}{n^2}\Big)\int_M|\nabla
R|^2\phi_r^2.
\end{align}
Using the Ricci identity, we have
\begin{align}\label{2-lemmas-9}
\mathring{R}_{kj,ij}\mathring{R}_{ki}=&(\mathring{R}_{kj,ji}+\mathring{R}_{pj}R_{pkij}
+\mathring{R}_{kp}R_{pjij})\mathring{R}_{ki}\notag\\
=&\frac{n-2}{2n}R_{,ki}\mathring{R}_{ki}+\mathring{R}_{pj}\mathring{R}_{ki}R_{pkij}
+\mathring{R}_{kp}\mathring{R}_{ki}\mathring{R}_{pi}+\frac{R}{n}|\mathring{\rm Ric}|^2\notag\\
=&\frac{n-2}{2n}R_{,ki}\mathring{R}_{ki}-W_{ijkl}\mathring{R}_{ik}\mathring{R}_{jl}
+\frac{n}{n-2}\mathring{R}_{ij}\mathring{R}_{jk}\mathring{R}_{ki}\notag\\
&+\frac{1}{n-1}R|\mathring{\rm Ric}|^2,
\end{align}
which shows
\begin{align}\label{2-lemmas-10}
-2a\int_M\mathring{R}_{kj,i}\mathring{R}_{ki,j}\phi_{r}^2=&2a\int_M\mathring{R}_{kj,ij}\mathring{R}_{ki}\phi_{r}^2
+2a\int_M\mathring{R}_{kj,i}\mathring{R}_{ki}(\phi_{r}^2)_{j}\notag\\
=&-2a\int_M\Big(W_{ijkl}\mathring{R}_{ik}\mathring{R}_{jl}-\frac{n}{n-2}\mathring{R}_{ij}\mathring{R}_{jk}
\mathring{R}_{ki}\notag\\
&-\frac{1}{n-1}R|\mathring{\rm Ric}|^2\Big)\phi_{r}^2+\frac{(n-2)a}{n}\int_MR_{,ki}\mathring{R}_{ki}\phi_{r}^2\notag\\
&+2a\int_M\mathring{R}_{kj,i}\mathring{R}_{ki}(\phi_{r}^2)_{j}.
\end{align}
Using the Cauchy inequality, we have
\begin{align}\label{2-lemmas-11}
2a\int_M\mathring{R}_{kj,i}\mathring{R}_{ki}(\phi_{r}^2)_{j}\leq\epsilon_1\int_M|\nabla \mathring{\rm Ric}|^2\phi_{r}^2+\frac{4a^2}{\epsilon_1}\int_M|\mathring{\rm Ric}|^2|\nabla \phi_{r}|^2
\end{align}
and
\begin{align}\label{2-lemmas-12}
&\frac{(n-2)a}{n}\int_MR_{,ki}\mathring{R}_{ki}\phi_{r}^2\notag\\
=&-\frac{(n-2)^2a}{2n^2}\int_M|\nabla R|^2\phi_{r}^2
-\frac{(n-2)a}{n}\int_MR_{,k}\mathring{R}_{ki}(\phi_{r}^2)_{i}\notag\\
\leq&\Big(\epsilon_2-\frac{(n-2)^2a}{2n^2}\Big)\int_M|\nabla R|^2\phi_{r}^2+\frac{(n-2)^2a^2}{n\epsilon_2}\int_M|\mathring{\rm Ric}|^2|\nabla \phi_{r}|^2,
\end{align}
where we used the inequality
\begin{align}\label{2-lemmas-13}
&-\frac{(n-2)a}{n}\int_MR_{,k}\mathring{R}_{ki}(\phi_{r}^2)_{i}\notag\\
=&-2\sum_{k}\int_M(R_{,k}\phi_{r})\Big(\frac{(n-2)a}{n}\sum_{i}\mathring{R}_{ki}(\phi_{r})_{i}\Big)\notag\\
\leq&\epsilon_2\sum_{k}\int_MR_{,k}^2\phi_{r}^2
+\frac{(n-2)^2a^2}{n\epsilon_2}\sum_{i,k}\int_M\Big(\mathring{R}_{ki}(\phi_{r})_{i}\Big)^2\notag\\
\leq&\epsilon_2\int_M|\nabla R|^2\phi_{r}^2+\frac{(n-2)^2a^2}{n\epsilon_2}\int_M|\mathring{\rm Ric}|^2|\nabla \phi_{r}|^2.
\end{align}
Substituting \eqref{2-lemmas-11} and \eqref{2-lemmas-12} into \eqref{2-lemmas-10} yields
\begin{align}\label{2-lemmas-14}
&-2a\int_M\mathring{R}_{kj,i}\mathring{R}_{ki,j}\phi_{r}^2\\
\leq&-2a\int_M\Big(W_{ijkl}\mathring{R}_{ik}\mathring{R}_{jl}-\frac{n}{n-2}\mathring{R}_{ij}\mathring{R}_{jk}
\mathring{R}_{ki}\notag\\
&-\frac{1}{n-1}R|\mathring{\rm Ric}|^2\Big)\phi_{r}^2+\epsilon_1\int_M|\nabla \mathring{\rm Ric}|^2\phi_{r}^2\notag\\
&+\Big(\epsilon_2-\frac{(n-2)^2a}{2n^2}\Big)\int_M|\nabla R|^2\phi_{r}^2+\Big(\frac{4}{\epsilon_1}\notag\\
&+\frac{(n-2)^2}{n\epsilon_2}\Big)a^2\int_M|\mathring{\rm Ric}|^2|\nabla \phi_{r}|^2,
\end{align}
which together with \eqref{2-lemmas-8} gives the desired estimate \eqref{2-lemmas-6}.
\endproof

For complete manifold $(M^n,g)$ with the metric $g$ critical for $\mathcal{F}_t$ over
$\mathscr{M}_1(M^n)$, we have the following result.

\begin{lem}\label{Lemma21}
Let $M^n$ be a complete manifold and $g$ be a critical for $\mathcal{F}_t$ over
$\mathscr{M}_1(M^n)$. Then for any two positive constants $\epsilon_3\in(0,1)$ and $\epsilon_4$, we have
\begin{align}\label{2-lemmas-15}
\int_M|\nabla \mathring{\rm Ric}|^2\phi_{r}^2\leq&\frac{1}{1-\epsilon_3}\int_M\Bigg[2W_{ijkl}\mathring{R}_{jl}\mathring{R}_{ik}
-\frac{4}{n-2}\mathring{R}_{ij}\mathring{R}_{jk}
\mathring{R}_{ki}\notag\\
&+\frac{2(n-2)+2n(n-1)t}{n(n-1)}R|\mathring{\rm Ric}|^2\notag\\
&+\Big(\epsilon_4+\frac{(n-2)(1+2t)}{2n}\Big)|\nabla R|^2\Bigg]\phi_{r}^2\notag\\
&+\frac{1}{1-\epsilon_3}\Big(\frac{n(1+2t)^2}{\epsilon_4}+\frac{1}{\epsilon_3}\Big)\int_M|\mathring{\rm Ric}|^2|\nabla \phi_{r}|^2.
\end{align}

\end{lem}

\proof It has been shown by
Catino in \cite[Proposition 2.1]{Catino2015} that a metric $g$ is critical for $\mathcal{F}_t$ over
$\mathscr{M}_1(M^n)$ if and only if it satisfies the following equations
\begin{align}
\Delta
\mathring{\rm R}_{ij}=&(1+2t)R_{,ij}-\frac{1+2t}{n} (\Delta R)
g_{ij}-2R_{ikjl}\mathring{\rm R}_{kl}\notag\\
&-\frac{2+2nt}{n}R \mathring{\rm R}_{ij} + \frac{2}{n} |\mathring{\rm Ric}|^{2} g_{ij},\label{2-lemmas-16}\\
[n+4(&n-1)t]\Delta R=(n-4)[| R_{ij}|^{2}+t
R^{2}-\lambda],\label{2-lemmas-17}
\end{align}
where $\lambda =\mathcal{F}_t (g)$.

It is easy to see from \eqref{2-lemmas-16} that
\begin{align}\label{2-lemmas-18}
\mathring{\rm R}_{ij}\Delta \mathring{\rm R}_{ij}=&(1+2t) \mathring{\rm R}_{ij}R_{,ij}-2R_{ikjl} \mathring{\rm R}_{kl}\mathring{\rm R}_{ij}
-\frac{2+2nt}{n} R|\mathring{\rm Ric}|^2\notag\\
=&(1+2t) \mathring{\rm R}_{ij}R_{,ij}-\frac{2(n-2)+2n(n-1)t}{n(n-1)} R|\mathring{\rm Ric}|^2\notag\\
&+\frac{4}{n-2}\mathring{\rm R}_{ij}\mathring{\rm R}_{jk}\mathring{\rm R}_{ki}-2W_{ikjl} \mathring{\rm R}_{kl}\mathring{\rm R}_{ij}.
\end{align}
Thus,
\begin{align}\label{2-lemmas-19}
\int_M|\nabla \mathring{\rm Ric}|^2\phi_{r}^2=&-\int_M\mathring{R}_{ij}\Delta\mathring{R}_{ij}\phi_{r}^2
-\int_M\mathring{R}_{ij}\mathring{R}_{ij,k}(\phi_{r}^2)_{k}\notag\\
=&\int_M\Big(2W_{ijkl}\mathring{R}_{jl}\mathring{R}_{ik}
-\frac{4}{n-2}\mathring{R}_{ij}\mathring{R}_{jk}
\mathring{R}_{ki}\notag\\
&+\frac{2(n-2)+2n(n-1)t}{n(n-1)}R|\mathring{\rm Ric}|^2\Big)\phi_{r}^2\notag\\
&-(1+2t)\int_M\mathring{\rm R}_{ij}R_{,ij}\phi_{r}^2
-\int_M\mathring{R}_{ij}\mathring{R}_{ij,k}(\phi_{r}^2)_{k}.
\end{align}
Applying the inequality
$$-\int_M\mathring{R}_{ij}\mathring{R}_{ij,k}(\phi_{r}^2)_{k}\leq\epsilon_3\int_M|\nabla \mathring{\rm Ric}|^2\phi_{r}^2+\frac{1}{\epsilon_3}\int_M|\mathring{\rm Ric}|^2|\nabla \phi_{r}|^2$$
and
$$\aligned
-(1+2t)\int_M\mathring{\rm R}_{ij}R_{,ij}\phi_{r}^2\leq&\Big(\epsilon_4+\frac{(n-2)(1+2t)}{2n}\Big)\int_M|\nabla R|^2\phi_{r}^2\\
&+\frac{n(1+2t)^2}{\epsilon_4}\int_M|\mathring{\rm Ric}|^2|\nabla \phi_{r}|^2
\endaligned$$
into \eqref{2-lemmas-19} gives the desired estimate \eqref{2-lemmas-15}.
\endproof

\begin{lem}\label{Lemma24}
On every Einstein manifold $(M^n,g)$, we have
\begin{align}\label{Weyl-1}
\frac{1}{2}\Delta |W|^2\geq&\frac{n+1}{n-1}|\nabla |W||^2+\frac{2}{n}R|W|^2-2C_n |W|^3,
\end{align}
where $C_n$ is defined by
\begin{equation}\label{cn}
C_n=\begin{cases}\frac{\sqrt{6}}{4},&\text{if\ }n=4;\\
\frac{4\sqrt{10}}{15},&\text{if\ }n=5;\\
\frac{n-2}{\sqrt{n(n-1)}} +\frac{n^2-n-4}{2\sqrt{(n-2)(n-1)n(n+1)}},&\text{if\ }n\geq 6.
\end{cases}
\end{equation}
In particular, if the scalar curvature of Einstein metric $g$ is positive,
then it is of constant positive sectional curvature, provided either
\begin{align}\label{Integ-1}
C_n|W|<\frac{1}{n}R,
\end{align}
or

(1) for $n\neq5$,
\begin{equation}\label{2-Pinchconstant}
\Big(\int_M|W|^{\frac{n}{2}}\Big)^{\frac{2}{n}}<E_n\,Q_g(M),
\end{equation}
where $E_n$ is given by \eqref{En};

(2) for $n=5$,
\begin{equation}\label{2-Pinchconstant-2}
\Big(\int_M|W|^{\frac{5}{2}}\Big)^{\frac{2}{5}}\leq\frac{2\sqrt{15}-4}{\sqrt{10}}\,Q_g(M).
\end{equation}

\end{lem}

\proof In \cite{Hebey1996} (or see \cite{FuXiao2017DGA}), it has been proved for an Einstein manifold that
\begin{align}\label{2-Sec-51}
\frac{1}{2}\Delta |W|^2=&|\nabla W|^2+\frac{2}{n}R|W|^2\notag\\
&-2(2W_{ijkl}W_{ipkq}W_{pjql}+\frac{1}{2}W_{ijkl}W_{klpq}W_{pqij})\notag\\
\geq&\frac{n+1}{n-1}|\nabla |W||^2+\frac{2}{n}R|W|^2\notag\\
&-2(2W_{ijkl}W_{ipkq}W_{pjql}+\frac{1}{2}W_{ijkl}W_{klpq}W_{pqij}),
\end{align}
where the inequality in \eqref{2-Sec-51}, we used the refined Kato inequality (see \cite{Calderbank2000}) of an Einstein manifold.

When $n=4,5$, we have (see \cite{FuXiao2017DGA,Catino2016})
\begin{align}\label{2-Sec-52}
|2W_{ijkl}W_{ipkq}W_{pjql}+\frac{1}{2}W_{ijkl}W_{klpq}W_{pqij}|\leq&\frac{\sqrt{6}}{4}|W|^3
\end{align}
and
\begin{align}\label{2-Sec-53}
|2W_{ijkl}W_{ipkq}W_{pjql}+\frac{1}{2}W_{ijkl}W_{klpq}W_{pqij}|\leq&\frac{4\sqrt{10}}{15}|W|^3,
\end{align}
respectively.

When $n=6$, making use of the inequality proved by Li and Zhao \cite{LZ1994} (see the formulas (13) and (14) in \cite{LZ1994}, or see \cite{CF2017}):
\begin{align}
2|W_{ijkl}W_{ipkq}W_{pjql}|=&|W_{ijkl}(W_{ipkq}W_{pjql}-W_{jpkq}W_{piql})|\notag\\
\leq&\frac{n-2}{\sqrt{n(n-1)}}|W|^3\notag
\end{align}
and Huisken \cite{Huisken1985}:
\begin{align}
|W_{ijkl}W_{klpq}W_{pqij}|\leq\frac{n^2-n-4}{\sqrt{(n-2)(n-1)n(n+1)}}|W|^3\notag
\end{align}
gives
\begin{align}\label{2-Sec-54}
|2W_{ijkl}&W_{ipkq}W_{pjql}+\frac{1}{2}W_{ijkl}W_{klpq}W_{pqij}|\notag\\
\leq&2|W_{ijkl}W_{ipkq}W_{pjql}|+\frac{1}{2}|W_{ijkl}W_{klpq}W_{pqij}|\notag\\
\leq&\frac{n-2}{\sqrt{n(n-1)}}+\frac{n^2-n-4}{2\sqrt{(n-2)(n-1)n(n+1)}}.
\end{align}
Hence, the desired estimate \eqref{Weyl-1} follows by inserting \eqref{2-Sec-52}-\eqref{2-Sec-54} into \eqref{2-Sec-51}.

Since the scalar curvature of Einstein metric $g$ is positive, $M^n$ must be of compact from the Myer's Theorem. The estimate \eqref{Integ-1} comes from integrating both sides of \eqref{Weyl-1}.

Let $v=|W|$. Then \eqref{Weyl-1} becomes
\begin{align}\label{2-Sec-55}
v\Delta v\geq&\frac{2}{n-1}|\nabla v|^2+\frac{2}{n}Rv^2-2C_n v^3,
\end{align}
and hence
\begin{align}\label{2-Sec-56}
v^{\alpha}\Delta v^{\alpha}=&v^{\alpha}[\alpha(\alpha-1)v^{\alpha-2}|\nabla v|^2+\alpha v^{\alpha-1}\Delta v]\notag\\
=&\Big(1-\frac{1}{\alpha}\Big)|\nabla v^{\alpha}|^2+\alpha v^{2\alpha-2}v\Delta v\notag\\
\geq&\Big(1-\frac{n-3}{(n-1)\alpha}\Big)|\nabla v^{\alpha}|^2-2\alpha C_nv^{2\alpha+1}+\frac{2\alpha}{n} Rv^{2\alpha}
\end{align}
which shows
\begin{align}\label{2-Sec-57}
0\geq&\Big(2-\frac{n-3}{(n-1)\alpha}\Big)\int_M|\nabla v^{\alpha}|^2-2\alpha C_n\int_Mv^{2\alpha+1}+\frac{2\alpha}{n} R\int_Mv^{2\alpha}.
\end{align}
Therefore, for $2-\frac{n-3}{(n-1)\alpha}>0$, by virtue of \eqref{1thm-Int-7} with $u=v^{\alpha}$, we have
\begin{align}\label{2-Sec-58}
0\geq&\Big(2-\frac{n-3}{(n-1)\alpha}\Big)Q_g(M)\Big(\int_Mv^{\frac{2n\alpha}{n-2}}\Big)^{\frac{n-2}{n}}-2\alpha C_n\int_Mv^{2\alpha+1}\notag\\
&+\Big[\frac{2\alpha}{n}-\frac{n-2}{4(n-1)}\Big(2-\frac{n-3}{(n-1)\alpha}\Big)\Big]R\int_Mv^{2\alpha}\notag\\
\geq&\Bigg[\Big(2-\frac{n-3}{(n-1)\alpha}\Big)Q_g(M)-2\alpha C_n\Big(\int_Mv^{\frac{n}{2}}\Big)^{\frac{2}{n}}\Bigg]\Big(\int_Mv^{\frac{2n\alpha}{n-2}}\Big)^{\frac{n-2}{n}}\notag\\
&+\Big[\frac{2\alpha}{n}-\frac{n-2}{4(n-1)}\Big(2-\frac{n-3}{(n-1)\alpha}\Big)\Big]R\int_Mv^{2\alpha}.
\end{align}

When $n\neq5$, taking $\frac{1}{\alpha}=\frac{n-1}{n-3}\Big(1+\sqrt{1-\frac{8(n-3)}{n(n-2)}}\Big)$ which satisfies
$$\frac{2\alpha}{n}-\frac{n-2}{4(n-1)}\Big(2-\frac{n-3}{(n-1)\alpha}\Big)=0.$$
Thus, if
\begin{align}\label{2-Sec-59}
\Big(\int_M|W|^{\frac{n}{2}}\Big)^{\frac{2}{n}}<&\frac{1}{2 C_n\alpha }\Big(2-\frac{n-3}{(n-1)\alpha}\Big)Q_g(M)\notag\\
=&\frac{4(n-1)}{n-2}\frac{Q_g(M)}{nC_n},
\end{align}
\eqref{2-Sec-58} shows that $W=0$ and it is of constant positive sectional curvature.
When $n=5$, \eqref{2-Sec-58} becomes
\begin{align}\label{2-Sec-60}
0\geq&\Bigg[\Big(2-\frac{1}{2\alpha}\Big)Q_g(M)-2\alpha C_5\Big(\int_Mv^{\frac{5}{2}}\Big)^{\frac{2}{5}}\Bigg]\Big(\int_Mv^{\frac{10\alpha}{3}}\Big)^{\frac{3}{5}}\notag\\
&+\Big[\frac{2\alpha}{5}-\frac{3}{16}\Big(2-\frac{1}{2\alpha}\Big)\Big]R\int_Mv^{2\alpha}.
\end{align}
Taking $\alpha=\frac{\sqrt{15}}{8}$ in \eqref{2-Sec-60} gives
\begin{align}\label{2-Sec-60-1}
0\geq&\Bigg[\frac{2\sqrt{15}-4}{\sqrt{15}}Q_g(M)-\frac{\sqrt{15}}{4} C_5\Big(\int_Mv^{\frac{5}{2}}\Big)^{\frac{2}{5}}\Bigg]\Big(\int_Mv^{\frac{10\alpha}{3}}\Big)^{\frac{3}{5}}\notag\\
&+\Big(\sqrt{\frac{3}{20}}-\frac{3}{8}\Big)R\int_Mv^{2\alpha},
\end{align}
which shows that $W=0$ and it is of constant positive sectional curvature provided \eqref{2-Pinchconstant-2} holds.

We complete the proof of Lemma \ref{Lemma24}. \endproof

The following lemma comes from \cite{Huang2018,FuXiao2015} (for the case of $\lambda=\frac{2}{n-2}$, see \cite{Huang2017,Catino2016, FuXiao2017Monath}):

\begin{lem}\label{Lemma23}
On every Riemannian manifold $(M^n,g)$, for any $\lambda\in \mathbb{R}$, the following estimate holds
\begin{align}\label{Maininequ-1}
&\Big|-W_{ijkl}\mathring{R}_{jl}\mathring{R}_{ik}
+\lambda\mathring{R}_{ij}\mathring{R}_{jk}\mathring{R}_{ki}\Big|\notag\\
\leq&
\sqrt{\frac{n-2}{2(n-1)}}\Big(|W|^2+\frac{2(n-2)\lambda^2}{n}|\mathring{\rm Ric}|^2
\Big)^{\frac{1}{2}}|\mathring{\rm Ric}|^2\notag\\
=&\sqrt{\frac{n-2}{2(n-1)}}\Big|W+\frac{\lambda}{\sqrt{2n}}\mathring{{\rm Ric}} \mathbin{\bigcirc\mkern-15mu\wedge} g\Big||\mathring{\rm Ric}|^2.
\end{align}

\end{lem}

\section{Proof of results}
\subsection{Proof of Theorem \ref{thm1-1}}
By combining \eqref{2-lemmas-6} with \eqref{2-lemmas-15}, we derive
\begin{align}\label{2-lemmas-20}
&\frac{1}{1-\epsilon_3}\int_M\Bigg[2W_{ijkl}\mathring{R}_{jl}\mathring{R}_{ik}
-\frac{4}{n-2}\mathring{R}_{ij}\mathring{R}_{jk}
\mathring{R}_{ki}\notag\\
&\hspace{.5cm}+\frac{2(n-2)+2n(n-1)t}{n(n-1)}R|\mathring{\rm Ric}|^2+\Big(\epsilon_4+\frac{(n-2)(1+2t)}{2n}\Big)|\nabla R|^2\Bigg]\phi_{r}^2\notag\\
&+\frac{1}{1-\epsilon_3}\Big(\frac{n(1+2t)^2}{\epsilon_4}+\frac{1}{\epsilon_3}\Big) \int_M|\mathring{\rm Ric}|^2|\nabla \phi_{r}|^2\notag\\
\geq&\frac{1}{a^2+1+\epsilon_1}\int_M\Bigg[2a\Big(
W_{ijkl}\mathring{R}_{jl}\mathring{R}_{ik}
-\frac{n}{n-2}\mathring{R}_{ij}\mathring{R}_{jk}
\mathring{R}_{ki}-\frac{1}{n-1}R|\mathring{\rm Ric}|^2\Big)\notag\\
&+|C^{(a,b)}|^2+\Big(-\frac{(n-2)^2(b-1)^2}{4n(n-1)^2}+\frac{(n-2)^2a}{2n(n-1)} -\epsilon_2\Big)|\nabla R|^2\Bigg]\phi_{r}^2\notag\\
&-\Big(\frac{4}{\epsilon_1}+\frac{(n-2)^2}{n\epsilon_2}\Big)\frac{a^2}{a^2+1+\epsilon_1} \int_M|\mathring{\rm Ric}|^2|\nabla \phi_{r}|^2,
\end{align}
which gives
\begin{align}\label{2-lemmas-21}
&\Big[\frac{1}{1-\epsilon_3}\Big(\frac{n(1+2t)^2}{\epsilon_4}+\frac{1}{\epsilon_3}\Big)
+\Big(\frac{4}{\epsilon_1}+\frac{(n-2)^2}{n\epsilon_2}\Big)\frac{a^2}{a^2+1+\epsilon_1}\Big]\int_M|\mathring{\rm Ric}|^2|\nabla \phi_{r}|^2\notag\\
\geq&-2\frac{a^2-(1-\epsilon_3)a+(1+\epsilon_1)}{(1-\epsilon_3)(a^2+1+\epsilon_1)}
\int_MW_{ijkl}\mathring{R}_{jl}\mathring{R}_{ik}\phi_{r}^2\notag\\
&+\frac{\frac{4}{n}a^2-2(1-\epsilon_3)a+\frac{4}{n}(1+\epsilon_1)}{(1-\epsilon_3)(a^2+1
+\epsilon_1)}\frac{n}{n-2}\int_M\mathring{R}_{ij}\mathring{R}_{jk}\mathring{R}_{ki}\phi_{r}^2\notag\\
&-\frac{2}{n(n-1)}\frac{na(1-\epsilon_3)+(a^2+1+\epsilon_1)[(n-2)+n(n-1)t]}{(1-\epsilon_3)(a^2+1
+\epsilon_1)}\int_MR|\mathring{\rm Ric}|^2\phi_{r}^2\notag\\
&+\Bigg[\frac{1}{a^2+1+\epsilon_1}\Big(-\frac{(n-2)^2(b-1)^2}{4n(n-1)^2}+\frac{(n-2)^2a}{2n(n-1)}-\epsilon_2\Big)\notag\\
&-\frac{1}{1-\epsilon_3}\Big(\epsilon_4+\frac{(n-2)(1+2t)}{2n}\Big)\Bigg]\int_M|\nabla R|^2\phi_{r}^2\notag\\
&+\frac{1}{a^2+1+\epsilon_1}\int_M|C^{(a,b)}|^2\phi_{r}^2.
\end{align}
In particular, for any positive $\epsilon_1,\epsilon_3$, we have $a^2-(1-\epsilon_3)a+(1+\epsilon_1)>0$. Hence, \eqref{2-lemmas-21} is equivalent to
\begin{align}\label{2-lemmas-22}
&\frac{(1-\epsilon_3)(a^2+1+\epsilon_1)}{a^2-(1-\epsilon_3)a+(1+\epsilon_1)}\Big[\frac{1}{1-\epsilon_3}
\Big(\frac{n(1+2t)^2}{\epsilon_4}+\frac{1}{\epsilon_3}\Big)\notag\\
&+\Big(\frac{4}{\epsilon_1}+\frac{(n-2)^2}{n\epsilon_2}\Big)\frac{a^2}{a^2+1+\epsilon_1}\Big]\int_M|\mathring{\rm Ric}|^2|\nabla \phi_{r}|^2\notag\\
\geq&-2\int_MW_{ijkl}\mathring{R}_{jl}\mathring{R}_{ik}\phi_{r}^2\notag\\
&+\frac{\frac{4}{n}a^2-2(1-\epsilon_3)a
+\frac{4}{n}(1+\epsilon_1)}{a^2-(1-\epsilon_3)a+(1+\epsilon_1)}\frac{n}{n-2}
\int_M\mathring{R}_{ij}\mathring{R}_{jk}\mathring{R}_{ki}\phi_{r}^2\notag\\
&-\frac{2}{n(n-1)}\frac{na(1-\epsilon_3)+(a^2+1+\epsilon_1)[(n-2)+n(n-1)t]}{a^2-(1-\epsilon_3)a
+(1+\epsilon_1)}\int_MR|\mathring{\rm Ric}|^2\phi_{r}^2\notag\\
&+\frac{(1-\epsilon_3)(a^2+1+\epsilon_1)}{a^2-(1-\epsilon_3)a+(1+\epsilon_1)} \Bigg[-\frac{1}{1-\epsilon_3}\Big(\epsilon_4+\frac{(n-2)(1+2t)}{2n}\Big)\notag\\
&\frac{1}{a^2
+1+\epsilon_1}\Big(-\frac{(n-2)^2(b-1)^2}{4n(n-1)^2}+\frac{(n-2)^2a}{2n(n-1)} -\epsilon_2\Big)\Bigg]\int_M|\nabla R|^2\phi_{r}^2\notag\\
&+\frac{1-\epsilon_3}{a^2-(1-\epsilon_3)a+(1+\epsilon_1)}\int_M|C^{(a,b)}|^2\phi_{r}^2.
\end{align}
Now, we fixe $\epsilon_1,\epsilon_3$ and minimize the function
$$f(a)=\frac{\frac{4}{n}a^2-2(1-\epsilon_3)a+\frac{4}{n}(1+\epsilon_1)}{a^2-(1-\epsilon_3)a+(1+\epsilon_1)}$$
by taking
$$a=\sqrt{1+\epsilon_1},$$
then \eqref{2-lemmas-22} becomes
\begin{align}\label{2-lemmas-23}
&\frac{\sqrt{1+\epsilon_1}(1-\epsilon_3)}{2\sqrt{1+\epsilon_1}-(1-\epsilon_3)}\Big[\frac{1}{1-\epsilon_3}
\Big(\frac{n(1+2t)^2}{\epsilon_4}+\frac{1}{\epsilon_3}\Big)\notag\\
&+\frac{1}{2}\Big(\frac{4}{\epsilon_1}+\frac{(n-2)^2}{n\epsilon_2}\Big)\Big]\int_M|\mathring{\rm Ric}|^2|\nabla \phi_{r}|^2\notag\\
\geq&-\int_MW_{ijkl}\mathring{R}_{jl}\mathring{R}_{ik}\phi_{r}^2
+\frac{\frac{4}{n}\sqrt{1+\epsilon_1}-(1-\epsilon_3)}{2\sqrt{1+\epsilon_1}-(1-\epsilon_3)}\frac{n}{n-2}
\int_M\mathring{R}_{ij}\mathring{R}_{jk}\mathring{R}_{ki}\phi_{r}^2\notag\\
&-\frac{1}{n(n-1)}\frac{n(1-\epsilon_3)+2\sqrt{1+\epsilon_1}[(n-2)+n(n-1)t]}{2\sqrt{1+
\epsilon_1}-(1-\epsilon_3)}\int_MR|\mathring{\rm Ric}|^2\phi_{r}^2\notag\\
&+\frac{\sqrt{1+\epsilon_1}(1-\epsilon_3)}{2\sqrt{1+\epsilon_1}-(1-\epsilon_3)}\Bigg[\frac{1}{2(1+\epsilon_1)}
\Big(-\frac{(n-2)^2(b-1)^2}{4n(n-1)^2}\notag\\
&+\frac{(n-2)^2}{2n(n-1)}\sqrt{1+\epsilon_1}-\epsilon_2\Big)-\frac{1}{1-\epsilon_3}\Big(\epsilon_4+\frac{(n-2)(1+2t)}{2n}\Big)\Bigg]\int_M|\nabla R|^2\phi_{r}^2\notag\\
&+\frac{1-\epsilon_3}{2\sqrt{1+\epsilon_1}[2\sqrt{1+\epsilon_1}-(1-\epsilon_3)]}\int_M|C^{(a,b)}|^2\phi_{r}^2.
\end{align}
Applying the inequality \eqref{Maininequ-1} with
$$\lambda=\frac{\frac{4}{n}\sqrt{1+\epsilon_1}-(1-\epsilon_3)}{2\sqrt{1+\epsilon_1}-(1-\epsilon_3)}\frac{n}{n-2}$$
into \eqref{2-lemmas-23} and taking $b=1$, we have
\begin{align}\label{2-lemmas-24}
&\frac{\sqrt{1+\epsilon_1}(1-\epsilon_3)}{2\sqrt{1+\epsilon_1}-(1-\epsilon_3)}\Big[\frac{1}{1-\epsilon_3}
\Big(\frac{n(1+2t)^2}{\epsilon_4}+\frac{1}{\epsilon_3}\Big)\notag\\
&+\frac{1}{2}\Big(\frac{4}{\epsilon_1}+\frac{(n-2)^2}{n\epsilon_2}\Big)\Big]\int_M|\mathring{\rm Ric}|^2|\nabla \phi_{r}|^2\notag\\
\geq&\int_M\Bigg[-\sqrt{\frac{n-2}{2(n-1)}}\Big|W+\frac{\sqrt{n}}{\sqrt{2}(n-2)}
\frac{\frac{4}{n}\sqrt{1+\epsilon_1}-(1-\epsilon_3)}{2\sqrt{1+\epsilon_1}-(1-\epsilon_3)}\mathring{{\rm Ric}} \mathbin{\bigcirc\mkern-15mu\wedge} g\Big|\notag\\
&-\frac{1}{n(n-1)}\frac{n(1-\epsilon_3)+2\sqrt{1+\epsilon_1}[(n-2)
+n(n-1)t]}{2\sqrt{1+\epsilon_1}-(1-\epsilon_3)}R\Bigg]|\mathring{\rm Ric}|^2\phi_{r}^2\notag\\
&+\frac{\sqrt{1+\epsilon_1}(1-\epsilon_3)}{2\sqrt{1+\epsilon_1}-(1-\epsilon_3)}\Bigg[\frac{1}{2(1+\epsilon_1)}
\Big(\frac{(n-2)^2}{2n(n-1)}\sqrt{1+\epsilon_1}-\epsilon_2\Big)\notag\\
&-\frac{1}{1-\epsilon_3}\Big(\epsilon_4+\frac{(n-2)(1+2t)}{2n}\Big)\Bigg]\int_M|\nabla R|^2\phi_{r}^2\notag\\
&+\frac{1-\epsilon_3}{2\sqrt{1+\epsilon_1}[2\sqrt{1+\epsilon_1}-(1-\epsilon_3)]}\int_M|C^{(a,b)}|^2\phi_{r}^2.
\end{align}

For all $\epsilon_1,\epsilon_3$, we have
\begin{equation}\label{2-lemmas-25}
\Big(\frac{\frac{4}{n}\sqrt{1+\epsilon_1}-(1-\epsilon_3)}{2\sqrt{1+\epsilon_1}-(1-\epsilon_3)}\Big)^2
=\Big(1-\frac{2-\frac{4}{n}}{2-\frac{1-\epsilon_3}{\sqrt{1+\epsilon_1}}}\Big)^2
<D_n,
\end{equation}
where $D_n$ is given by \eqref{dn}.
Using the fact that $W$ is perpendicular to $\mathring{{\rm Ric}} \mathbin{\bigcirc\mkern-15mu\wedge} g$, we have
\begin{align}\label{2-lemmas-26}
\Big|W&+\frac{\sqrt{n}}{\sqrt{2}(n-2)}
\frac{\frac{4}{n}\sqrt{1+\epsilon_1}-(1-\epsilon_3)}{2\sqrt{1+\epsilon_1}-(1-\epsilon_3)}\mathring{{\rm Ric}} \mathbin{\bigcirc\mkern-15mu\wedge} g\Big|^2\notag\\
&=|W|^2+\frac{n}{2(n-2)^2}\Big(\frac{\frac{4}{n}\sqrt{1+\epsilon_1}-(1-\epsilon_3)}{2\sqrt{1+\epsilon_1}-(1-\epsilon_3)}\Big)^2|\mathring{{\rm Ric}} \mathbin{\bigcirc\mkern-15mu\wedge} g|^2\notag\\
&<|W|^2+\frac{n D_n}{2(n-2)^2}|\mathring{{\rm Ric}} \mathbin{\bigcirc\mkern-15mu\wedge} g|^2\notag\\
&=\Big|W+\frac{\sqrt{n D_n}}{\sqrt{2}(n-2)}\mathring{{\rm Ric}} \mathbin{\bigcirc\mkern-15mu\wedge} g\Big|^2,
\end{align}
which shows
\begin{align}\label{2-lemmas-27}
\Big|W&+\frac{\sqrt{n}}{\sqrt{2}(n-2)}
\frac{\frac{4}{n}\sqrt{1+\epsilon_1}-(1-\epsilon_3)}{2\sqrt{1+\epsilon_1}-(1-\epsilon_3)}\mathring{{\rm Ric}} \mathbin{\bigcirc\mkern-15mu\wedge} g\Big|\notag\\
&<\Big|W+\frac{\sqrt{n D_n}}{\sqrt{2}(n-2)}\mathring{{\rm Ric}} \mathbin{\bigcirc\mkern-15mu\wedge} g\Big|.
\end{align}
On the other hand, if
\begin{align}\label{2-lemmas-29}
n+2[(n-2)+n(n-1)t]\leq0,
\end{align}
then we have
\begin{align}\label{2-lemmas-28}
2[(n-2)+n(n-1)t]<&\frac{n(1-\epsilon_3)+2\sqrt{1+\epsilon_1}[(n-2)
+n(n-1)t]}{2\sqrt{1+\epsilon_1}-(1-\epsilon_3)}\notag\\
&=\frac{n\frac{1-\epsilon_3}{\sqrt{1+\epsilon_1}}+2[(n-2)
+n(n-1)t]}{2-\frac{1-\epsilon_3}{\sqrt{1+\epsilon_1}}}\notag\\
&<\frac{n+2[(n-2)
+n(n-1)t]}{2-\frac{1-\epsilon_3}{\sqrt{1+\epsilon_1}}}\notag\\
&<\frac{n+2[(n-2)
+n(n-1)t]}{2}
\end{align}
and hence
\begin{align}\label{2-lemmas-30}
\frac{n(1-\epsilon_3)+2\sqrt{1+\epsilon_1}[(n-2)
+n(n-1)t]}{2\sqrt{1+\epsilon_1}-(1-\epsilon_3)}<0.
\end{align}
Moreover, taking $\epsilon_2=\epsilon_1\sqrt{1+\epsilon_1}$ and $\epsilon_4=\frac{\epsilon_1(1-\epsilon_3)}{2\sqrt{1+\epsilon_1}}$, we have
\begin{align}\label{1-add-2-lemmas-31}
\frac{(n-2)^2}{4n(n-1)}&\frac{1-\epsilon_3}{\sqrt{1+\epsilon_1}}-\epsilon_2\frac{1-\epsilon_3}{2(1+\epsilon_1)}
-\epsilon_4-\frac{(n-2)(1+2t)}{2n}\notag\\
=&\Big(\frac{(n-2)^2}{4n(n-1)}-\epsilon_1\Big)\frac{1-\epsilon_3}{\sqrt{1+\epsilon_1}}
-\frac{(n-2)(1+2t)}{2n}.
\end{align}
Since,
\begin{align}\label{2-add-2-lemmas-31}
\Big(\frac{(n-2)^2}{4n(n-1)}-\epsilon_1\Big)\frac{1-\epsilon_3}{\sqrt{1+\epsilon_1}}\rightarrow \frac{(n-2)^2}{4n(n-1)}
\end{align}
as $\epsilon_1\rightarrow 0$ and $\epsilon_3\rightarrow 0$, there exist two positive constants $\overline{\epsilon}_1,\overline{\epsilon}_3$ such that
\begin{align}\label{2-lemmas-31}
\frac{(n-2)^2}{4n(n-1)}&\frac{1-\overline{\epsilon}_3}{\sqrt{1+\overline{\epsilon}_1}}-\overline{\epsilon}_2\frac{1-\overline{\epsilon}_3}{2(1+\overline{\epsilon}_1)}
-\overline{\epsilon}_4-\frac{(n-2)(1+2t)}{2n}\notag\\
>&\frac{6(n-2)^2}{25n(n-1)}-\frac{(n-2)(1+2t)}{2n}.
\end{align}
Hence, if
\begin{align}\label{2-lemmas-32}
\frac{6(n-2)^2}{25n(n-1)}-\frac{(n-2)(1+2t)}{2n}\geq0
\end{align}
then \eqref{2-lemmas-31} shows
\begin{align}\label{3-add2-lemmas-31}
\frac{(n-2)^2}{4n(n-1)}&\frac{1-\overline{\epsilon}_3}{\sqrt{1+\overline{\epsilon}_1}}-\overline{\epsilon}_2\frac{1-\overline{\epsilon}_3}{2(1+\overline{\epsilon}_1)}
-\overline{\epsilon}_4-\frac{(n-2)(1+2t)}{2n}>0
\end{align}
which is equivalent to
\begin{align}\label{2-lemmas-33}
\frac{1}{2(1+\overline{\epsilon}_1)}&
\Big(\frac{(n-2)^2}{2n(n-1)}\sqrt{1+\overline{\epsilon}_1}-\overline{\epsilon}_2\Big)\notag\\
&-\frac{1}{1-\overline{\epsilon}_3}\Big(\overline{\epsilon}_4+\frac{(n-2)(1+2t)}{2n}\Big)>0.
\end{align}
It is easy to check that the scope of $t$ satisfying both \eqref{2-lemmas-29} and \eqref{2-lemmas-32} is equivalent to \eqref{1111thm-Int-1}. Therefore, if \eqref{1111thm-Int-1} and \eqref{1thm-Int-4} both hold, then we have
\begin{align}\label{2-lemmas-34}
&\frac{\sqrt{1+\overline{\epsilon}_1}(1-\overline{\epsilon}_3)}{2\sqrt{1+\overline{\epsilon}_1}-(1-\overline{\epsilon}_3)}
\Big[\frac{1}{1-\overline{\epsilon}_3}
\Big(n(1+2t)^2\frac{2\sqrt{1+\overline{\epsilon}_1}}{\overline{\epsilon}_1(1-\overline{\epsilon}_3)}
+\frac{1}{\overline{\epsilon}_3}\Big)\notag\\
&+\frac{1}{2}\Big(\frac{4}{\overline{\epsilon}_1}+\frac{(n-2)^2}{n\overline{\epsilon}_1\sqrt{1+\overline{\epsilon}_1}}
\Big)\Big]\int_M|\mathring{\rm Ric}|^2|\nabla \phi_{r}|^2\notag\\
\geq&\int_M\Bigg[-\sqrt{\frac{n-2}{2(n-1)}}\Big|W+\frac{\sqrt{n}}{\sqrt{2}(n-2)}
\frac{\frac{4}{n}\sqrt{1+\overline{\epsilon}_1}-(1-\overline{\epsilon}_3)}{2\sqrt{1+\overline{\epsilon}_1}
-(1-\overline{\epsilon}_3)}\mathring{{\rm Ric}} \mathbin{\bigcirc\mkern-15mu\wedge} g\Big|\notag\\
&-\frac{1}{n(n-1)}\frac{n(1-\overline{\epsilon}_3)+2\sqrt{1+\overline{\epsilon}_1}[(n-2)
+n(n-1)t]}{2\sqrt{1+\overline{\epsilon}_1}-(1-\overline{\epsilon}_3)}R\Bigg]|\mathring{\rm Ric}|^2\phi_{r}^2\notag\\
\geq&0
\end{align}
from \eqref{2-lemmas-24}.

Since
\begin{equation}\label{3-Sec-35}
\int_M|\mathring{\rm Ric}|^2<\infty,
\end{equation}
then we have
\begin{equation}\label{3-Sec-36}\aligned
\int_M|\mathring{\rm Ric}|^2|\nabla \phi_{r}|^2\rightarrow0,
\endaligned\end{equation}
as $r\rightarrow \infty$, which together with \eqref{2-lemmas-34} shows that $M^n$ is Einstein. In this case, \eqref{1thm-Int-4}
becomes
\begin{align}\label{3-Sec-37}
|W|&\leq-\frac{2\sqrt{2}[(n-2)+n(n-1)t]}{n\sqrt{(n-1)(n-2)}}R,
\end{align}
which gives
\begin{align}\label{3-Sec-38}
C_n|W|&\leq-\frac{2\sqrt{2}[(n-2)+n(n-1)t]}{n\sqrt{(n-1)(n-2)}}C_nR.
\end{align}
When $n=3$, we have $W=0$ automatically. When $n\geq4$ and $t$ satisfies \eqref{1111thm-Int-6}, we have
\begin{align}\label{3-Sec-39}
-\frac{2\sqrt{2}[(n-2)+n(n-1)t]}{n\sqrt{(n-1)(n-2)}}C_nR&<\frac{1}{n}R
\end{align}
which, combining with \eqref{Integ-1}, shows that $W=0$ and hence
$M^n$ is of constant positive sectional curvature.

\subsection{Proof of Theorem \ref{thm1-2}}

From \eqref{Maininequ-1}, it is easy to see
\begin{align}\label{4-sec-1}
2W_{ijkl}&\mathring{R}_{jl}\mathring{R}_{ik}
-\frac{4}{n-2}\mathring{R}_{ij}\mathring{R}_{jk}
\mathring{R}_{ki}\notag\\
\leq&\sqrt{\frac{2(n-2)}{n-1}}\Big|W+\frac{\sqrt{2}}{\sqrt{n}(n-2)}\mathring{{\rm Ric}} \mathbin{\bigcirc\mkern-15mu\wedge} g\Big||\mathring{\rm Ric}|^2.
\end{align}
Applying \eqref{4-sec-1} into \eqref{2-lemmas-15} and using the Kato inequality, we obtain
\begin{align}\label{4-sec-2}
\int_M|\nabla |\mathring{\rm Ric}||^2\phi_{r}^2
\leq&\int_M|\nabla \mathring{\rm Ric}|^2\phi_{r}^2\notag\\
\leq&\frac{1}{1-\epsilon_3}\int_M\Bigg[\sqrt{\frac{2(n-2)}{n-1}}\Big|W+\frac{\sqrt{2}}{\sqrt{n}(n-2)}\mathring{{\rm Ric}} \mathbin{\bigcirc\mkern-15mu\wedge} g\Big|\notag\\
&+\frac{2(n-2)+2n(n-1)t}{n(n-1)}R\Bigg]|\mathring{\rm Ric}|^2\phi_{r}^2\notag\\
&+\Big(\epsilon_4+\frac{(n-2)(1+2t)}{2n}\Big)\int_M|\nabla R|^2\phi_{r}^2\notag\\
&+\frac{1}{1-\epsilon_3}\Big(\frac{n(1+2t)^2}{\epsilon_4}+\frac{1}{\epsilon_3}\Big)\int_M|\mathring{\rm Ric}|^2|\nabla \phi_{r}|^2.
\end{align}
Taking $u=|\mathring{\rm Ric}|\phi_{r}$ in \eqref{1thm-Int-7} and applying \eqref{4-sec-2} yield
\begin{align}\label{4-sec-3}
Q_g(M)&\Big(\int_M(|\mathring{\rm Ric}|\phi_{r})^{\frac{2n}{n-2}}\Big)^{\frac{n-2}{n}}\notag\\
\leq&\int_M\Big(|\nabla (|\mathring{\rm Ric}|\phi_{r})|^2+\frac{n-2}{4(n-1)}R|\mathring{\rm Ric}|^2\phi_{r}^2\Big)\notag\\
\leq&(1+\epsilon_5)\int_M|\nabla |\mathring{\rm Ric}||^2\phi_{r}^2+\Big(1+\frac{1}{\epsilon_5}\Big)\int_M|\mathring{\rm Ric}|^2|\nabla \phi_{r}|^2\notag\\
&+\frac{n-2}{4(n-1)}\int_MR|\mathring{\rm Ric}|^2\phi_{r}^2\notag\\
\leq&\frac{1+\epsilon_5}{1-\epsilon_3}\sqrt{\frac{2(n-2)}{n-1}}\int_M\Big|W+\frac{\sqrt{2}}{\sqrt{n}(n-2)}\mathring{{\rm Ric}} \mathbin{\bigcirc\mkern-15mu\wedge} g\Big||\mathring{\rm Ric}|^2\phi_{r}^2\notag\\
&+\frac{1}{n-1}\Big[\frac{n-2}{4}+\frac{2(n-2)+2n(n-1)t}{n}\frac{1+\epsilon_5}{1-\epsilon_3}\Big]\int_MR|\mathring{\rm Ric}|^2\phi_{r}^2\notag\\
&+(1+\epsilon_5)\Big[\frac{1}{\epsilon_5}+\frac{1}{1-\epsilon_3}\Big(\frac{n(1+2t)^2}{\epsilon_4}+\frac{1}{\epsilon_3}\Big)\Big]\int_M|\mathring{\rm Ric}|^2|\nabla \phi_{r}|^2\notag\\
&+(1+\epsilon_5)\Big(\epsilon_4+\frac{(n-2)(1+2t)}{2n}\Big)\int_M|\nabla R|^2\phi_{r}^2.
\end{align}
Inserting the following H\"{o}lder inequality
$$\aligned
\int_M&\Big|W+\frac{\sqrt{2}}{\sqrt{n}(n-2)}\mathring{{\rm Ric}} \mathbin{\bigcirc\mkern-15mu\wedge} g\Big||\mathring{\rm Ric}|^2\phi_{r}^2\\
\leq&\Big(\int_M\Big|W+\frac{\sqrt{2}}{\sqrt{n}(n-2)}\mathring{{\rm Ric}} \mathbin{\bigcirc\mkern-15mu\wedge} g\Big|^{\frac{n}{2}}\Big)^{\frac{2}{n}}\Big(\int_M(|\mathring{\rm Ric}|\phi_{r})^{\frac{2n}{n-2}}\Big)^{\frac{n-2}{n}}
\endaligned$$
into \eqref{4-sec-3} deduces
\begin{align}\label{4-sec-4}
&\Bigg[Q_g(M)-\frac{1+\epsilon_5}{1-\epsilon_3}\sqrt{\frac{2(n-2)}{n-1}} \Big(\int_M\Big|W+\frac{\sqrt{2}}{\sqrt{n}(n-2)}\mathring{{\rm Ric}} \mathbin{\bigcirc\mkern-15mu\wedge} g\Big|^{\frac{n}{2}}\Big)^{\frac{2}{n}}\Bigg]\times\notag\\
&\hspace{8cm}\times\Big(\int_M(|\mathring{\rm Ric}|\phi_{r})^{\frac{2n}{n-2}}\Big)^{\frac{n-2}{n}}\notag\\
\leq&\frac{1}{n-1}\Big[\frac{n-2}{4}+\frac{2(n-2)+2n(n-1)t}{n}\frac{1+\epsilon_5}{1-\epsilon_3}\Big]
\int_MR|\mathring{\rm Ric}|^2\phi_{r}^2\notag\\
&+(1+\epsilon_5)\Big[\frac{1}{\epsilon_5}+\frac{1}{1-\epsilon_3}\Big(\frac{n(1+2t)^2}{\epsilon_4}
+\frac{1}{\epsilon_3}\Big)\Big]\int_M|\mathring{\rm Ric}|^2|\nabla \phi_{r}|^2\notag\\
&+(1+\epsilon_5)\Big(\epsilon_4+\frac{(n-2)(1+2t)}{2n}\Big)\int_M|\nabla R|^2\phi_{r}^2.
\end{align}
We check that if $t<-\frac{1}{2}$, then $(n-2)+n(n-1)t<0$ and
$$\frac{-8[(n-2)+n(n-1)t]}{n(n-2)}>1.$$
Hence for all $\epsilon_3,\epsilon_5$, it holds that
\begin{align}\label{4-sec-5}
\frac{1-\epsilon_3}{1+\epsilon_5}<\frac{-8[(n-2)+n(n-1)t]}{n(n-2)},
\end{align}
which is equivalent to
\begin{align}\label{4-sec-6}
\frac{n-2}{4}+\frac{2(n-2)+2n(n-1)t}{n}\frac{1+\epsilon_5}{1-\epsilon_3}<0.
\end{align}
If \eqref{222thm-Int-3} occurs, then
\begin{align}\label{4-sec-7}
\frac{\Big(\int_M\Big|W+\frac{\sqrt{2}}{\sqrt{n}(n-2)}\mathring{{\rm Ric}} \mathbin{\bigcirc\mkern-15mu\wedge} g\Big|^{\frac{n}{2}}\Big)^{\frac{2}{n}}}{\sqrt{\frac{n-1}{2(n-2)}}Q_g(M)}<1.
\end{align}
In this case, there exist $\tilde{\epsilon}_3,\tilde{\epsilon}_5$ such that
\begin{align}\label{4-sec-8}
\frac{1-\tilde{\epsilon}_3}{1+\tilde{\epsilon}_5} =\frac{\Big(\int_M\Big|W+\frac{\sqrt{2}}{\sqrt{n}(n-2)}\mathring{{\rm Ric}} \mathbin{\bigcirc\mkern-15mu\wedge} g\Big|^{\frac{n}{2}}\Big)^{\frac{2}{n}}}{\sqrt{\frac{n-1}{2(n-2)}}Q_g(M)},
\end{align}
which, from \eqref{4-sec-4}, shows that
\begin{align}\label{4-sec-9}
0\leq&\frac{1}{n-1}\Big[\frac{n-2}{4}+\frac{2(n-2)+2n(n-1)t}{n}\frac{1+\tilde{\epsilon}_5} {1-\tilde{\epsilon}_3}\Big]\int_MR|\mathring{\rm Ric}|^2\phi_{r}^2\notag\\
&+(1+\tilde{\epsilon}_5)\Big[\frac{1}{\tilde{\epsilon}_5}+\frac{1} {1-\tilde{\epsilon}_3}\Big(\frac{2n^2(1+2t)^2}{-(n-2)(1+2t)}
+\frac{1}{\tilde{\epsilon}_3}\Big)\Big]\times\notag\\
&\hspace{6cm}\times\int_M|\mathring{\rm Ric}|^2|\nabla \phi_{r}|^2
\end{align}
by taking $\epsilon_4=-\frac{(n-2)(1+2t)}{2n}$. It follows from \eqref{222thm-Int-1} that
\begin{align}\label{4-sec-10}
\int_M|\mathring{\rm Ric}|^2|\nabla \phi_{r}|^2\rightarrow 0
\end{align}
as $r\rightarrow \infty$, which together with \eqref{4-sec-9} gives that $M^n$ is Einstein.

Since $M^n$ is Einstein, then \eqref{222thm-Int-3} becomes
\begin{equation}\label{4-sec-11}
\Big(\int_M|W|^{\frac{n}{2}}\Big)^{\frac{2}{n}}<\sqrt{\frac{n-1}{2(n-2)}}Q_g(M).
\end{equation}
When $n=4,5$, we can check that
\begin{equation}\label{4-sec-12}
\sqrt{\frac{n-1}{2(n-2)}}<E_n,
\end{equation}
which implies from Lemma \ref{Lemma24} that $W=0$ and hence $M^n$ is of constant positive sectional curvature.
When $n\geq 6$, it is easy to check
\begin{equation}\label{4-sec-14}
\frac{4}{n}\sqrt{\frac{2(n-1)}{n-2}}<\frac{n-2}{\sqrt{n(n-1)}}+\frac{n^2-n-4}{2\sqrt{(n-2)(n-1)n(n+1)}}
\end{equation}
which is equivalent to
\begin{equation}\label{4-sec-13}
E_n<\sqrt{\frac{n-1}{2(n-2)}}.
\end{equation}
Therefore, if \eqref{222thm-Int-4} holds, then $M^n$ is also of constant positive sectional curvature.

We complete the proof of Theorem \ref{thm1-2}.

\bibliographystyle{Plain}

\end{document}